\documentclass[11pt]{amsart}
\usepackage{amsmath}
  \usepackage{url,graphicx}
  \usepackage{amssymb, color, pstricks-add, manfnt}
  \usepackage{amsmath, amsthm,  amsfonts}
  \usepackage[plainpages=false]{hyperref}

  \theoremstyle{plain}
  \newtheorem{theorem}{Theorem}[section]
  \newtheorem{lemma}{Lemma}[section]
  \newtheorem{remark}{Remark}[section]

  \numberwithin{equation}{section}
  \numberwithin{figure}{section}

\renewcommand{\baselinestretch}{1.00}
\begin{document}

\title{On the Dirichlet problem for Monge-Amp\`ere type equations}

\author{Feida Jiang}
\address{School of Science, Nanjing University of Science $\&$
Technology, Nanjing 210094, P.R.China} \email{jfd2001@163.com}

\author{Neil S. Trudinger}
\address{Centre for Mathematics and Its Applications, The Australian National University,
              Canberra ACT 0200, Australia}
\email{Neil.Trudinger@anu.edu.au}

\author{Xiao-Ping Yang}
\address{School of Science, Nanjing University of Science $\&$
Technology, Nanjing 210094, P.R.China}
\email{yangxp@mail.njust.edu.cn}

\thanks{This research was finished when the first author was visiting the
Australian National University in 2011. It was supported by National
Natural Science Foundation of China(No.11071119) and the Australian
Research Council.}


\date{\today}


\maketitle

\abstract {In this paper, we prove second derivative estimates
together with classical solvability for the Dirichlet problem of
certain Monge-Amp\`ere type equations under sharp hypotheses. In
particular we assume that the matrix function in the augmented
Hessian is regular in the sense used by Trudinger and Wang in their
study of global regularity in optimal transportation
\cite{TruWang2009} as well as the existence of a smooth subsolution.
The latter hypothesis replaces a barrier condition also used in
their work. The applications to optimal transportation and
prescribed Jacobian equations are also indicated.}
\endabstract


\baselineskip=12.8pt
\parskip=3pt
\renewcommand{\baselinestretch}{1.38}

\section{Introduction}

\vskip10pt

This paper is concerned with existence of globally smooth solutions
to the Dirichlet problem for certain Monge-Amp\`ere type equations.
These equations arise in many applications, notably in optimal
transportation, geometric optics and conformal geometry. The
Dirichlet problem of the Monge-Amp\`ere type equations under
consideration has the following general form
\begin{equation}\label{Dirichlet BVP of OTE}
\left\{ \begin{array}{rcll}
    \det\{D^2u-A(x,Du)\} & = & B(x,Du),     &  x \in \Omega ,           \\
                       u & = & \varphi(x),  &  x \in \partial\Omega ,   \\
         \end{array}  \right.
\end{equation}
where $\Omega$ is a bounded domain with smooth boundary in $n$
dimensional Euclidean space $\mathbb{R}^n$, $\varphi$ is a smooth
function on $\partial \Omega$. $A$ is a given $n \times n$ symmetric
matrix function defined on $\Omega \times \mathbb{R}^n$, and $B$ is
a scalar valued function defined on $\Omega \times \mathbb{R}^n$. We
use $(x,p)$ to define points in $\Omega \times \mathbb{R}^n$ so that
$A(x,p) \in \mathbb{R}^n\times \mathbb{R}^n$, $B(x,p) \in \mathbb{R}$. 
Also
$Du$ and $D^2u$ denote the gradient vector and Hessian matrix
respectively of a function $u\in C^2(\Omega)$ . The equation in
(\ref{Dirichlet BVP of OTE}) reduces to the standard Monge-Amp\`ere
equation when $A \equiv 0$. A solution $u\in C^2(\Omega)$ of
(\ref{Dirichlet BVP of OTE}) is called an elliptic solution whenever
$D^2u-A(x,Du)>0$, which implies $B>0$.

These equations have attracted significant interest in recent years;
for a recent survey and the earlier history see, for example,
\cite{TruWang2008}. The Heinz-Lewy example in \cite{Schulz90} shows
that there is no $C^1$ regularity for the equation in
(\ref{Dirichlet BVP of OTE}) without restrictions on $A$. When
$A\equiv0$, the existence of globally smooth solutions to the
Dirichlet problem (\ref{Dirichlet BVP of OTE}) for smooth uniformly
convex domains was finally settled by Caffarelli, Nirenberg and
Spruck \cite{CNSI84,CNS87}, Krylov \cite{Krylov83} and Ivochkina
\cite{Ivochkina80,Ivochkina83}. For application to optimal
transportation with quadratic cost function, the natural boundary
condition is the prescription of the gradient image $\Omega^* =
Du(\Omega)$ which is equivalent to an oblique boundary condition for
elliptic solutions. In this case local regularity for convex
$\Omega^*$ was proved by Caffarelli \cite{Caffarelli92JAMS} and the
corresponding global regularity for smooth uniformly convex domains
$\Omega$ and target domains $\Omega^*$ proved by Delano\"e
\cite{D91}, Caffarelli \cite{Caffarelli96Annals} and Urbas
\cite{Urbas97}. The first breakthrough for general optimal
transportation problems was obtained by Ma, Trudinger and Wang
\cite{MTW2005} who obtained an {\it a priori} interior second order
estimate under an analytical structure condition on the matrix $A$,
which they  called A3, together with local regularity results under
a generalized target convexity condition. These were subsequently
extended to global estimates and regularity in \cite{TruWang2009}
under a weak form of this condition called A3w; (see \cite{Tru2006}
for further discussion). Following this Loeper \cite{Loeper09Acta}
found a geometric characterization enabling him to adapt the proof
of necessity of target convexity  in \cite{MTW2005} to prove that
A3w is necessary for regularity in optimal transportation.

In this paper we show that the condition A3w also suffices for
global regularity in the Dirichlet problem not just in the optimal
transportation case but also for general  equations of the form
(\ref{Dirichlet BVP of OTE}). In this generality it is convenient to
adopt the terminology from \cite{Tru2006} and call the matrix
function $A$ regular if
\begin{equation}\label{A3W}
A_{ij,kl}(x,p)\xi_i\xi_j\eta_k\eta_l\geqslant 0,
\end{equation}
for all $(x,p)\in\Omega\times\mathbb{R}^n$,
$\xi,\eta\in\mathbb{R}^n$, $\xi\perp\eta$, where
$A_{ij,kl}=D^2_{p_kp_l}A_{ij}$.

We now  formulate the main results of this paper. The first result
is a global bound for second derivatives of equation (\ref{Dirichlet
BVP of OTE}) which replaces the barrier condition in
\cite{TruWang2009} by the existence of a subsolution.

\begin{theorem}\label{Th1}
Let $u\in C^4(\Omega)\cap C^2(\bar\Omega)$ be an elliptic solution
of problem (\ref{Dirichlet BVP of OTE}) in $\Omega$, where $B\in
C^2(\bar\Omega\times\mathbb{R}^n)$, $\inf B>0$ and $A\in
C^2(\bar\Omega\times\mathbb{R}^n)$ is regular. Suppose also there
exists a subsolution $\underline u\in C^2(\bar\Omega)$ of equation
(\ref{Dirichlet BVP of OTE}). Then we have the estimate
\begin{equation}\label{interior bound}
\sup\limits_\Omega |D^2u| \leqslant C(1+\sup\limits_{\partial
\Omega} |D^2u|),
\end{equation}
where the constant $C$ depends on $n, A, B, \Omega, \underline u$
and $\sup\limits_\Omega{(|u|+|Du|)}$.
\end{theorem}

From Theorem \ref{Th1}, we can infer a global bound for solutions of
the Dirichlet problem (\ref{Dirichlet BVP of OTE}). For this we
adopt the approach from \cite{Guanbo93,Guanbo95} and instead of
explicit geometric assumptions on the boundary $\partial\Omega$ as
formulated for example in \cite{Tru2006}, we assume more generally
the existence of a subsolution with the given boundary trace.

\begin{theorem}\label{Th2}
In addition to the assumptions in Theorem \ref{Th1}, suppose the
subsolution $\underline u = \varphi$ on $\partial\Omega$ with
$\varphi \in C^4(\bar\Omega), \partial\Omega\in C^4$. Then any
elliptic solution $u\in C^4(\Omega)\cap C^2(\bar\Omega))$ of the
Dirichlet problem (\ref{Dirichlet BVP of OTE}) satisfies the global
a priori estimate
\begin{equation}
\sup\limits_\Omega |D^2u|\leqslant C,
\end{equation}
where the constant $C$ depends on $n, A, B, \Omega, \underline u$
and $\sup\limits_\Omega{(|u|+|Du|)}$.
\end{theorem}

Existence theorems follow from Theorem \ref{Th2} by standard
methods. In general we need supplementary conditions on the matrix
function $A$ to control gradients and upper bounds of solutions.
These can be reduced in the optimal transportation case, where the
matrix $A$ is generated by a {\it cost function}, namely a function
$c\in C^4(\mathbb{R}^n\times\mathbb{R}^n)$ satisfying the
conditions:

\noindent (A1) \hspace{3mm}    For each $x\in\Omega$, $p\in
\mathbb{R}^n$, there exists a unique $y = Y(x,p)$ such that $D_xc(x,y)=p$;\\

\noindent (A2) \hspace{3mm}    $\det D^2_{x,y}c \neq 0$,
for all $x\in\Omega, y = Y(x,p), p\in \mathbb{R}^n$. \\

The matrix $A$ is then given by
\begin{equation}\label{OT}
A(x,p) = D_x^2c(x,Y(x,p)).
\end{equation}

We formulate here a particular theorem for the general case which
also embraces the main examples in optimal transportation. For this
we assume a structure condition:
\begin{equation}\label{S1}
A(x,p) \geqslant - \mu_0(1 + |p|^2) I,
\end{equation}
\noindent for all $x\in\Omega, p\in \mathbb{R}^n$ and some positive
constant $\mu_0$ and that the maximum eigenvalue of $A(x,0)$ is
non-negative. The last hypothesis can be replaced by a weaker
condition that there exists a bounded viscosity supersolution to
equation (\ref{Dirichlet BVP of OTE}).

We then have the following existence theorem.

\begin{theorem}\label{Th3}
Under the assumptions in Theorem \ref{Th2} together with the above
structure conditions,
 there exists a unique  elliptic classical solution $u\in C^3(\bar \Omega)$ of the Dirichlet problem
(\ref{Dirichlet BVP of OTE}) with $u \geqslant \underline u$.
\end{theorem}

The paper is organized as follows: In Section 2 we prove a
fundamental lemma about barriers for the linearized operator and use
it to establish the global second derivative bound in Theorem
\ref{Th1}. We also prove an analogous variant of the interior
Pogorelev estimate of Liu and Trudinger in \cite{Liu10pog}. In
Section 3, we derive the boundary estimates needed to complete the
proof of Theorem \ref{Th2}. Our treatment of these boundary
estimates
 also
provides the proof of the corresponding results  stated in
\cite{Tru2006} under the barrier condition introduced in
\cite{TruWang2009}.  In Section 4, we establish the necessary
gradient and solution bounds for our existence theorems. Once we
have the relevant {\it a priori} estimates, Theorem 1.3 follows
immediately from the method of continuity in [11, Chapter 17].
Finally, we discuss separately the special cases of optimal
transportation and prescribed Jacobian equations.

\section{Pogorelov type estimates}\label{Section 2}

In this section we prove global and interior Pogorelov type
estimates for second derivatives of equation (1.1). In particular we
use the regularity of the matrix function $A$ to construct a barrier
function for the linearized operator from a subsolution. The
resultant second derivative estimates then follow readily from the
corresponding proofs in \cite{TruWang2009} and \cite{Liu10pog}.

Let $F[u]=\log (\det (D^2u-A(x,Du)))$, it is known that $F$ is a
concave operator with respect to $D^2u-A(x,Du)$. The linearized
operator of $F$ is defined by
\begin{equation}\label{linearized operator L}
L=F^{ij}(D_{ij}-D_{p_k}A_{ij}(x,Du)D_k),
\end{equation}
where $F^{ij}=\frac{\partial F}{\partial w_{ij}}$ and
$\{F^{ij}\}=\{w^{ij}\}$ denotes the inverse matrix of
$\{w_{ij}\}\triangleq\{u_{ij}-A_{ij}(x,Du)\}$. Assume the function
$A_{ij}(x,p)\in C^2(\bar\Omega\times\mathbb{R}^n)$,
$i,j=1,\cdots,n$, $B(x,p)\in C^2(\bar\Omega\times\mathbb{R}^n)$,
$\tilde B(x,p)=\log(B(x,p))$ and set
\begin{equation}\label{linearized operator}
\mathcal {L}=L-\tilde B_{p_i}D_i.
\end{equation}
We introduce a fundamental lemma here:
\begin{lemma}\label{Lemma}
Under the above assumptions, suppose $u$ is an elliptic solution of
(\ref{Dirichlet BVP of OTE}), $\underline u$ is a strict elliptic
subsolution of equation (\ref{Dirichlet BVP of OTE}), if $A$ is
regular, then
\begin{equation}\label{key inequality}
\mathcal {L} \left({e^{K(\underline u-u)}}\right) \geqslant
\epsilon_1\sum\limits_{i} F^{ii}-C,
\end{equation}
holds in $\Omega$ for positive constant $K$ sufficiently large and
uniform positive constants $\epsilon_1, C$.
\end{lemma}
\noindent \textbf{Proof}. Since $\underline u$ is a strict
subsolution of (\ref{Dirichlet BVP of OTE}), $\underline u$
satisfies
\begin{displaymath}
F[\underline u]=\log (\det (D^2\underline u-A(x,D\underline
u))\geqslant \log(B(x,D\underline u)+\delta_0),
\end{displaymath}
for some $\delta_0>0$.

For any $x_0\in \Omega$, let $\underline u_\epsilon=\underline
u-\frac{\epsilon}{2}|x-x_0|^2$. For $\epsilon$ small enough, the
perturbation function $\underline u_\epsilon$ is still a strict
subsolution and satisfies
\begin{displaymath}
F[\underline u_\epsilon]=\log (\det (D^2\underline
u_\epsilon-A(x,D\underline u_\epsilon)))\geqslant
\log(B(x,D\underline u_\epsilon)+\tau),
\end{displaymath}
for some positive constant $\tau$.

Let $v=\underline u-u$, $v_\epsilon=\underline u_\epsilon-u$. By
calculation, we have
\begin{equation}\label{Lv}
\begin{array}{lll}
Lv&=&L(v_\epsilon)+L(\frac{\epsilon}{2}|x-x_0|^2)\\
  &=&\epsilon F^{ii}-\epsilon
  F^{ij}D_{p_k}A_{ij}(x,Du)(x-x_0)_k+F^{ij}\{D_{ij}v_\epsilon-D_{p_k}A_{ij}(x,Du)D_kv_\epsilon\}\\
  &=&\epsilon F^{ii}-\epsilon
  F^{ij}D_{p_k}A_{ij}(x,Du)(x-x_0)_k+
  F^{ij}\{D_{ij}(\underline u_\epsilon -u)-[A_{ij}(x,D\underline
  u_\epsilon)-A_{ij}(x,Du)]\}\\
  && +F^{ij}\{A_{ij}(x,D\underline
  u_\epsilon)-A_{ij}(x,Du)-D_{p_k}A_{ij}(x,Du)D_kv_\epsilon\}.
\end{array}
\end{equation}
By the concavity of $F$, we have
\begin{displaymath}
F[\underline u_\epsilon]-F[u]\leqslant F^{ij}\{D_{ij}(\underline
u_\epsilon -u)-[A_{ij}(x,D\underline u_\epsilon)-A_{ij}(x,Du)]\}.
\end{displaymath}
By the Taylor expansion, for some $\theta\in (0,1)$, we have
\begin{displaymath}
\begin{array}{ll}
&A_{ij}(x,D\underline
u_\epsilon)-A_{ij}(x,Du)-D_{p_k}A_{ij}(x,Du)D_kv_\epsilon\\
=&D_{p_k}A_{ij}(x,\hat{p})D_kv_\epsilon-D_{p_k}A_{ij}(x,Du)D_kv_\epsilon\\
=&\frac{\theta}{2}A_{ij,kl}(x,\bar p)D_kv_\epsilon
D_lv_\epsilon,\\
\end{array}
\end{displaymath}
where $\hat{p}=(1-\theta)Du+\theta D\underline u_\epsilon$, $\bar
p=(1-\bar \theta)Du+\bar \theta D\underline u_\epsilon$ and $\bar
\theta\in (0, \theta)$.

Thus, at $x=x_0$, we have
\begin{displaymath}
\begin{array}{lll}
Lv & \geqslant & \epsilon F^{ii}+F[\underline
u_\epsilon]-F[u]+\frac{\theta}{2}F^{ij}A_{ij,kl}(x,\bar p)D_kv D_lv\\
& = & \epsilon F^{ii}+\frac{\theta}{2}F^{ij}A_{ij,kl}(x,\bar
p)D_kv D_lv+ \log (B(x,D\underline u)+\tau)-\log (B(x,Du))\\
& \geqslant & \epsilon F^{ii}+\frac{\theta}{2}F^{ij}A_{ij,kl}(x,\bar
p)D_kv D_lv -C_1,
\end{array}
\end{displaymath}
where $C_1$ is a constant depends on $B$, $Du$, and $D\underline u$.

Let $\phi=e^{Kv}$ with positive constant $K$ to be determined, we
have
\begin{displaymath}
\begin{array}{lll}
L\phi & = & F^{ij}(D_{ij}\phi-D_{p_k}A_{ij}(x,Du)D_k\phi)\\
      & = & Ke^{Kv}Lv + K^2e^{Kv}F^{ij}D_ivD_jv\\
      & \geqslant & Ke^{Kv}\{\epsilon
F^{ii}+\frac{\theta}{2}F^{ij}A_{ij,kl}(x,\bar p)D_kv D_lv
-C_1+KF^{ij}D_ivD_jv\}.
\end{array}
\end{displaymath}
Without loss of generality, we assume that $Dv=(D_1v, 0, \cdots,
0)$. Since $A$ is regular, we have
\begin{displaymath}
\begin{array}{lll}
L\phi & \geqslant & Ke^{Kv}\{\epsilon
F^{ii}+\frac{\theta}{2}F^{ij}A_{ij,11}(x,\bar p)(D_1v)^2
+KF^{11}(D_1v)^2-C_1\}\\
& \geqslant & Ke^{Kv}\{\epsilon
F^{ii}+\frac{\theta}{2}\sum\limits_{i \hspace{1mm} or \hspace{1mm} j
=1}F^{ij}A_{ij,11}(x,\bar p)(D_1v)^2
+KF^{11}(D_1v)^2-C_1\}.\\
\end{array}
\end{displaymath}
Since the matrix $\{F^{ij}\}$ is positive definite, any $2\times 2$
diagonal minor has positive determinant. This implies
\begin{displaymath}
|F^{1i}|^2\leqslant F^{11}F^{ii}.
\end{displaymath}
Then the Cauchy inequality leads to the following inequality
\begin{displaymath}
|F^{1i}|\leqslant \eta F^{ii}+\frac{1}{4\eta}F^{11},
\end{displaymath}
for positive constant $\eta$.

Thus, we have
\begin{displaymath}
\begin{array}{lll}
L\phi & \geqslant & Ke^{Kv}\{\epsilon F^{ii}-\frac{\theta}{2}\eta
F^{ii}|A_{1i,11}(x,\bar
p)|(D_1v)^2-\frac{\theta}{8\eta}F^{11}|A_{1i,11}(x,\bar p)|(D_1v)^2
+KF^{11}(D_1v)^2-C_1\}.\\
\end{array}
\end{displaymath}
Choosing $\eta$ small such that $\eta \leqslant
\frac{\epsilon}{\theta \max{\{|A_{1i,11}(x,\bar p)|(D_1v)^2\}}}$ and
$K$ large such that $K\geqslant \frac{\theta \max|A_{1i,11}(x,\bar
p)|}{8\eta}$, we obtain
\begin{displaymath}
\begin{array}{lll}
L\phi & \geqslant & Ke^{Kv}\{
\frac{\epsilon}{2}F^{ii}-C_1\}.\\
\end{array}
\end{displaymath}
Thus, we have
\begin{displaymath}
\mathcal {L}\phi=L\phi-\tilde B_{p_i}D_i\phi\geqslant Ke^{Kv}\{
\frac{\epsilon}{2}F^{ii}-C_1\}-\tilde B_{p_i}D_i\phi,
\end{displaymath}
which leads to the conclusion of Lemma \ref{Lemma} with
$\epsilon_1=\min_{\bar \Omega}\{\frac{\epsilon}{2}K e^{Kv}\}$ and
$C=\max_{\bar \Omega}\{C_1Ke^{Kv}+\tilde B_{p_i}D_i\phi\}$.

\qed

\begin{remark}\label{remark2.1}
The function $\phi={e^{K(\underline u-u)}}$ is global barrier
function for the linearized operator $\mathcal {L}$. The inequality
$\mathcal {L} \left({e^{K(\underline u-u)}}\right) \geqslant
\epsilon_1\sum\limits_{i} F^{ii}-C$ in $\Omega$ will be used to
control the second order global and interior estimates in this
section and the boundary estimates in the next section.

\end{remark}

\begin{remark}
It is a fairly standard  calculation that a non-strict classical
subsolution of a uniformly elliptic partial differential equation
can be made strict using the linearized operator and the mean value
theorem.
 For example, if $\underline u$ is a non-strict subsolution,
then $\underline u + a e^{bx_1}$ is a strict subsolution for small
constant $a$ and large constant $b$, (see [11], Chapter 3).
 This can also be done in a neighbourhood
of the boundary, preserving boundary conditions by replacing $x_1$
by a defining function or more specifically the distance function.
By using the mean value theorem we do not need concavity or
convexity and the small constant $a$ controls the uniform
ellipticity near the subsolution. Hence we need only assume the
existence of a non-strict subsolution in Lemma 2.1; the inequality
(\ref{key inequality}) will still hold for the corresponding strict
subsolution. Thus, the second order {\it{apriori}} estimates will
also hold under the existence of a non-strict subsolution and we
only need to assume a non-strict subsolution in the hypotheses for
our theorems.
\end{remark}

\begin{remark}
In the previous paper \cite{TruWang2009}, to obtain the second
derivative bounds, a kind of global barrier condition, called {\it
A-boundedness} in \cite{Tru2006} is assumed, namely $\Omega$ is
called $A$-bounded with respect to $u$ if there exists a function
$\varphi\in C^2(\bar \Omega)$ satisfying
\begin{equation}\label{2.4}
[D_{ij}\varphi -
D_{p_k}A_{ij}(\cdot,Du)D_k\varphi]\xi_i\xi_j\geqslant |\xi|^2,
\end{equation}
in $\Omega$, for all $\xi\in \mathbb{R}^n$. From Lemma 2.1 it
 follows that $A$-boundedness can be replaced by the existence of
a subsolution for second derivative estimates.
\end{remark}

\begin{remark}
If $F$ is substituted by another concave operator
$F[u]=(\det(D^2u-A(x,Du)))^{\frac 1n}$, the conclusion of this lemma
still holds. The only difference is the form of
$F^{ij}=\frac{\partial F}{\partial w_{ij}}$. We have $F^{ij}=
\frac{1}{n} f w^{ij}$ if $F[u]=(\det (D^2u-A(x,Du)))^{\frac{1}{n}}$,
while $F^{ij}=w^{ij}$ if $F[u]=\log(D^2u-A(x,Du))$. Here we make a
convention that $\{w^{ij}\}$ denotes the inverse matrix of
$\{w_{ij}\}\triangleq\{u_{ij}-A_{ij}(x,Du)\}$.
\end{remark}

Now from Lemma 2.1, we obtain the global Pogorelov estimate, Theorem
1.1, by appropriate adjustment to the proof of Theorem 3.1 in
\cite{TruWang2009}. This is done by substituting the barrier
function $\widetilde{\varphi} = \epsilon_1^{-1}\exp(K(\underline u -
u))$ in place of the barrier from condition (\ref{2.4}). That is we
maximize the function,
\begin{equation}\label{glp}
\exp\{\frac{a}{2}|Du|^2+b\phi\}w_{\xi\xi},
\end{equation}
over $\Omega$ and $|\xi| = 1$,where
$w_{\xi\xi}=w_{ij}\xi_i\xi_j=(u_{ij}-A_{ij})\xi_i\xi_j$,
$\phi=e^{K(\underline u-u)}$ with $K$ as in Lemma 2.1, and $a$, $b$
are positive constants to be determined. Then we arrive again at the
estimates (3.11) and (3.13) in \cite{TruWang2009}, but now with an
additional dependence on $\underline u$, and Theorem 1.1 follows.

As remarked in the introduction we may similarly modify the proof of
the interior Pogorelev estimate in  \cite{Liu10pog} to obtain the
following variant.

\begin{theorem}\label{Th2.1}
Let $u\in C^4(\Omega)\cap C^{0,1}(\bar\Omega)$ be an elliptic
solution of problem (\ref{Dirichlet BVP of OTE}) in $\Omega$, where
$B\in C^2(\bar\Omega\times\mathbb{R}^n)$, $\inf B>0$ and $A\in
C^2(\bar\Omega\times\mathbb{R}^n)$ is regular. Suppose also
 there exists an elliptic  subsolution $\underline u\in C^2(\bar\Omega)$
of equation (\ref{Dirichlet BVP of OTE}) and a degenerate elliptic
supersolution $u_0\in C^{1,1}(\Omega)\cap C^{0,1}(\bar\Omega)$ such
that $u = u_0$ on $\partial\Omega$. Then we have for any
$\Omega'\Subset\Omega$
\begin{equation}\label{intest}
\sup_{\Omega'}|D^2u|\leqslant C,
\end{equation}
where the constant $C$ depends on $n,A,B,\Omega,\Omega', \underline
u, u_0$,and $\sup_\Omega(|u|+|Du|)$.
\end{theorem}

Note that if $A(x,0) = 0$, then our estimate (\ref{intest}) agrees
with the usual Pogorelov estimate \cite{GT2001}, by taking $u_0 =
0$.

\section{ Boundary estimates}\label{Section 3}

In order to complete the proof of Theorem \ref{Th2}, we need to
obtain an {\it a priori} estimate for the Hessian $D^2u$ on the
boundary $\partial \Omega$ for solutions of the Dirichlet problem
(\ref{Dirichlet BVP of OTE}). First, we note from Remark 2.2 that we
can assume the subsolution $\underline u$ is strict provided we
restrict to a neighbourhood of $\partial \Omega$. As remarked in the
previous section, this can be shown by adding a function
$c_1\exp(c_2d(x))$, where $d(x)=dist(x,\partial \Omega)$ denotes the
distance function in $\Omega$, and $c_1$ and $c_2$ are positive
constants. Accordingly, we can retain Lemma \ref{Lemma} in a
neighbourhood of $\partial \Omega$, which will suffice for boundary
estimates.

Next, we need to observe that the regularity of the matrix function
$A$ will be preserved under coordinate changes, that is, under a
coordinate change, equation (\ref{Dirichlet BVP of OTE}) is
transformed into an equation of the same form with transformed
matrix $A$ satisfying (\ref{A3W}) with respect to the transformed
gradient variables. Specifically if we take a diffeomorphism, $y =
\psi(x)$ with Jacobian matrix $J = \{\psi_{ij}\} = \{D_j\psi_i\}$,
then
$$ w_{ij}(x) = \psi_{ki}\psi_{lj}D_{y_ky_l}u + D_i(\psi_{kj})D_{y_k}u
 - A_{ij}(x, JD_yu).$$
The transformed matrix function $A^\prime$ and scalar $B^\prime$ are
thus given by
\begin{equation}
\left\{ \begin{array}{rcll}
A{^\prime}_{ij}(\cdot,p)& = &\psi^{ik} \psi^{jl} [A_{kl}(\psi^{-1}, Jp) - D_k(\psi_{sl})p_s],\\
B{^\prime}(\cdot,p)& = & (\det J)^{-2} B(\psi^{-1}, Jp), \\
\end{array}  \right.
\end{equation}
where $\{\psi^{ij}\} = J^{-1}$. The invariance of condition (1.2)
follows readily from (3.1).

Consequently, we may fix a point $x_0\in
\partial\Omega$, which we take as the origin, and a neighbourhood
$\mathcal {N}$ of $x_0$, such that $T=\mathcal{N}\cap
\partial \Omega$ lies in the hyperplane $\{x_n = 0\}$,
and the positive $x_n$ axis points into $\Omega$.
 We then have by differentiation,
$D_{\alpha\beta}u(x',0)=D_{\alpha\beta}\underline u(x',0)$,
$\alpha,\beta=1,\cdots,n-1$, $x'=(x_1,\cdots,x_{n-1})$, which leads
to the double tangential derivative estimate
$|D_{\alpha\beta}u(x)|\leqslant C$, $\alpha,\beta=1,\cdots,n-1$, on
$T$.

We then estimate the mixed tangential-normal derivatives $D_{\alpha
n}u$, $\alpha=1,\cdots,n-1$, on $T$ by differentiating the equation
(\ref{Dirichlet BVP of OTE}) and using the barrier in Remark 2.1.
With this barrier in hand, the derivation of these estimates  is
similar to the special case of the standard Monge-Amp\`ere equation,
as treated for example in [11]. For completeness and later
reference, we carry out the details here. First by
 differentiation of the  equation
\begin{equation}
F[u]=\log (\det (D^2u-A(x,Du)))=\tilde B(x,Du),
\end{equation}
we have
\begin{equation}
F^{ij}(u_{\alpha i j}-D_\alpha
A_{ij}(x,Du)-D_{p_k}A_{ij}(x,Du)u_{\alpha k})=\tilde
B_\alpha(x,Du)+\tilde B_{p_k}(x,Du)u_{\alpha k}, \ \ \alpha=1,
\cdots, n,
\end{equation}
which leads to
\begin{equation}
\mathcal {L}D_\alpha u =\tilde B_\alpha (x,Du)+F^{ij}D_\alpha
A_{ij}(x,Du),
\end{equation}
where $\mathcal {L}$ is defined by (\ref{linearized operator}). Then
we have for tangential derivatives,

\begin{equation}
    |\mathcal{L}(D_\alpha (u-\underline u))|\leqslant C(1+\sum_i F^{ii}), \ \ \alpha=1,\cdots,n-1,\ \  {\rm in} \ \Omega.
\end{equation}
We also observe that, on $\partial \Omega$ close to the origin, the
tangential derivative satisfies $D_\alpha(u-\underline u)=0$,
$\alpha=1,2,\cdots,n-1$, i.e., $|D_\alpha(u-\underline u)|\leqslant
C|x'|^2$. This implies
\begin{displaymath}
|D_\alpha(u-\underline u)|\leqslant C|x|^2, \ \ {\rm on} \ \
\partial\Omega_\delta,
\end{displaymath}
where $\Omega_\delta=\Omega \cap B_\delta (0)$ is a small
neighborhood of the origin, $\partial \Omega_\delta=(\partial
\Omega\cap B_\delta(0))\cup (\Omega \cap \partial B_\delta(0))$. We
always choose $\delta$ so small that $\partial \Omega \cap
B_\delta(0) \subset T$.

Let $v=1-\phi$, where $\phi=e^{K(\underline u-u)}$ is the barrier
function in Remark 2.1, then we have
\begin{equation}
    \mathcal {L}v\leqslant -\epsilon_1\sum\limits_{i} F^{ii}+C, \ \   {\rm in}\ \  \Omega,\ \ \  {\rm and}\ \ \
    v= 0, \ \  {\rm on}\ \  \partial \Omega,
\end{equation}
for positive constants $\epsilon_1$ and $C$. We consider a function
of the form
\begin{equation}
\psi = v + \mu x_n - k{x_n}^2,
\end{equation}
where $\mu$ and $k$ are positive constants to be determined. We then
have
\begin{equation}\label{phi}
    \mathcal {L}\psi \leqslant  -\frac{\epsilon_1}{4}(1+\sum\limits_i F^{ii}),\ \ {\rm in}\ \  \Omega_\delta,\ \ \ {\rm and}\ \ \ \psi \geqslant 0,  \ \ {\rm on}\ \  \partial \Omega_\delta,
\end{equation}
for $k$ sufficiently large and $\mu$, $\delta$ sufficiently small.

By modification of the function $\psi$, we employ a new barrier
function
\begin{equation}
\tilde \psi =a \psi + b|x|^2,
\end{equation}
with positive $a$ and $b$ to be determined. By a direct computation,
for $a \gg b \gg 1$, we have
\begin{equation}
    \mathcal{L}\tilde \psi =a\mathcal{L}\psi + b\mathcal{L}|x|^2\leqslant -(\frac{a\epsilon_1}{4}-b)(1+\sum_i F^{ii}),
\end{equation}
\begin{equation}
    |D_\alpha(u-\underline u)|=0\leqslant \tilde \psi,\ \ {\rm on} \ \partial \Omega\cap B_\delta(0),
\end{equation}
and
\begin{equation}
    |D_\alpha(u-\underline u)|\leqslant b\delta^2\leqslant a\psi+b\delta^2=\tilde \psi,\ \ {\rm on} \ \Omega\cap \partial B_\delta(0).
\end{equation}
Therefore, for $a\gg b\gg 1$ and $\delta \ll 1$, there holds
\begin{equation}
\begin{array}{rllll}
    |\mathcal{L}D_\alpha(u-\underline u)|+\mathcal{L}\tilde \psi & \leqslant & 0,           & {\rm in} & \Omega_\delta,\\
    |D_\alpha(u-\underline u)|              & \leqslant & \tilde \psi, & {\rm on} & \partial \Omega_\delta.
\end{array}
\end{equation}
By the maximum principle, we have
\begin{equation}
    |D_\alpha(u-\underline u)| \leqslant \tilde \psi,\ {\rm in}  \ \Omega_\delta.
\end{equation}
At the origin, we have $|D_\alpha(u-\underline u)(0)|=\tilde
\psi(0)=0$. For $|D_\alpha(u-\underline u)|$, taking $x'=0$,
dividing by $x_n$ and letting $x_n\rightarrow 0$, we get
\begin{equation}
    |D_{\alpha n}(u-\underline u)(0)|\leqslant a \psi_n(0)\leqslant C.
\end{equation}
Thus, we obtain the mixed tangential-normal derivative estimate
$|D_{\alpha n}u(x)|\leqslant C$, $\alpha=1,\cdots,n-1$, on $T$.

It therefore remains to estimate the double normal derivative
$D_{nn}u$ and for this the regularity of the matrix function $A$ is
critical. Following the idea in \cite{Tru95Acta} as used in
\cite{Tru95Notes}, we fix a unit vector $\xi\in \mathbb{R}^{n-1}$
and set, (with respect to our transformed coordinates),
\begin{equation}\label{3.1}
\begin{array}{lll}
w:&=&w[u]=w_{\alpha \beta}\xi_\alpha\xi_\beta\\
 &=&[D_{\alpha\beta}u-A_{\alpha\beta}(x,Du)]\xi_\alpha\xi_\beta\\
 &=&[D_{\alpha\beta}\varphi
 -A_{\alpha\beta}(x,D^\prime\varphi,D_nu)]\xi_\alpha\xi_\beta \  {\rm on}\ T,\\
\end{array}
\end{equation}
with $w > 0$ by virtue of the ellipticity of $u$. Here $D^\prime =
(D_1,......D_{n-1})$ denotes the tangential gradient. An estimate
from above for $D_{nn}u$ on $T$ will follow from equation
(\ref{Dirichlet BVP of OTE}), provided we obtain a positive lower
bound for $w$ on $T$. To get this we note that for sufficiently
large $K$ the function $\tilde w =w+K|x|^2$ will take a minimum on
$T\times S^{n-1}$ at some point $\bar x\in T$, and for some unit
vector $\bar \xi$. Extending $\varphi$ by defining $\varphi(x',x_n)
= \varphi(x',0)$ for $x_n >0$, it follows from the regularity of
$A$,  that the corresponding extension of $w$ is a concave function
of $D_nu$. Consequently from the differentiated equation (3.4), for
$\alpha=n$, we have
\begin{equation}
\mathcal {L}\tilde w \leqslant C(1+\sum_i F^{ii}).
\end{equation}
Now using the barrier in Remark \ref{remark2.1} as above, we obtain
this time a one-sided estimate $D_nw(\bar x)\geqslant -C$, that is
\begin{equation}\label{3.2}
D_{nn}u(\bar x)D_{p_n}A_{\alpha\beta}(\bar x,D^\prime \varphi,D_n
u)\bar \xi_\alpha\bar \xi_\beta \leqslant C.
\end{equation}
To reach our desired estimate for $D_{nn}u$, we then need to get a
positive lower bound for the coefficient in (\ref{3.2}). By the
ellipticity of the subsolution $\underline u$, we can fix a positive
constant $\delta$ for which
\begin{equation}\label{3.3}
w[\underline u]\geqslant \delta,\hspace{2mm}{\rm
for}\hspace{1.5mm}{\rm all}\hspace{1.5mm}x\in T, |\xi|=1.
\end{equation}
We also have
\begin{equation}\label{3.4}
\begin{array}{lll}
w[u]-w[\underline u]&=&- [A_{\alpha\beta}(\bar
x,D^\prime\varphi,D_nu)-A_{\alpha\beta}(\bar
x,D^\prime\varphi,D_n\underline
u)]\bar \xi_\alpha \bar \xi_\beta \\
&\geqslant & -D_n(u-\underline u)(\bar x)D_{p_n}A_{\alpha\beta}(\bar
x,D^\prime\varphi,D_nu)\bar \xi_\alpha \bar \xi_\beta,
\end{array}
\end{equation}
since $D_nu > D_n\underline u$ from the strictness of $\underline u$
and  again using the regularity of $A$. We also have an upper bound
\begin{equation}\label{3.5}
D_n(u-\underline u)\leqslant \kappa,
\end{equation}
for a positive constant $\kappa$.

Combining (\ref{3.3}), (\ref{3.4}) and (\ref{3.5}), we thus obtain
for $w[u]<\frac{\delta}{2}$,
\begin{equation}\label{3.6}
D_{p_n}A_{\alpha\beta}(\bar
x,D^\prime\varphi,D_nu)\bar\xi_\alpha\bar\xi_\beta>\frac{\delta}{2\kappa}.
\end{equation}
Hence, we conclude from (\ref{3.2}) an estimate $D_{nn}u(\bar
x)\leqslant C$ for a further constant $C$. Now utilizing equation
(\ref{Dirichlet BVP of OTE}) again, we obtain an estimate from below
for $w(\bar x)$ and finally an estimate from above for
$D_{nn}u(x_0)$. Since an estimate from below automatically follows
from the ellipticity of $u$, we complete the estimation of $D^2u$ on
$\partial \Omega$ and the proof of Theorem \ref{Th2}.

\begin{remark}
In the above proof, we could have used, in place of (\ref{3.1}), the
function
\begin{equation}\label{3.7}
w[u]=\left\{\det [w_{ij}]\right\}^{\frac{1}{n-1}}
\end{equation}
where indices $i,j=1,\cdots,n-1$, in direct accordance with the
technique in \cite{Tru95Acta}.
\end{remark}

\begin{remark}
In \cite{Tru2006}, the concept of domain $A$-convexity is introduced
extending that of $c$-convexity in optimal transportation. In
particular for $A$ a given $n \times n$ symmetric matrix function
defined on $\Omega \times \mathbb{R}^n$ we may define $\Omega$ to be
uniformly $A$-convex, with respect to $u$  if
\begin{equation}
[D_i \gamma_j + D_{p_k}A_{ij}(x,Du) \gamma_k] \tau_i \tau_j
\geqslant \delta_0
\end{equation}
for all $x \in \partial\Omega$, unit outer normal $\gamma$ and unit
tangent vector $\tau$ and some positive constant $\delta_0$. It then
follows
 that $\Omega$ is uniformly $A$-convex if
and only if, for any constant $K > 0$, there exists a defining
function $\varphi\in C^2(\bar \Omega)$ satisfying $\varphi = 0$ on
$\partial\Omega$ together with the inequality (\ref{2.4}) in a
neighbourhood of $\partial\Omega$. By using $\varphi$ in place of
the barrier in Remark \ref{remark2.1}, we conclude an estimate for
the Hessian $D^2u$ on $\partial\Omega$ for solutions of  the
Dirichlet problem \ref{Dirichlet BVP of OTE} for uniformly
$A$-convex domains, thereby proving Theorem 2.2 in \cite{Tru2006}.
 Note that the concept of $A$-convexity is also invariant
with respect to coordinate changes. Note that the replacement of
uniform $A$-convexity by the existence of a strict subsolution with
the same boundary trace for second derivative estimates, as
formulated in Theorem 2.1, is also pointed out in \cite{Tru2006}.
\end{remark}

\begin{remark}
Theorem \ref{Th2}, provides bounds for both $D^2u$ and the augmented
matrix function $w$. By the equation (\ref{Dirichlet BVP of OTE}),
we have the positive lower bounds for $w$ in $\bar\Omega$. Thus, the
uniform ellipticity of the operator $F$, with respect to $u$,
follows easily from the positive upper and lower bounds of $w$.
\end{remark}

\section{Existence theorem and some applications}\label{Section 4}

In this section, we complete the proof of  the classical solvability
result for the Dirichlet problem (\ref{Dirichlet BVP of OTE}) and
consider the applications to the optimal transportation and
prescribed Jacobian equations.

First we establish the necessary solution bounds and gradient bounds
for Theorem 1.3. By the comparison principle, we have $u\geqslant
\underline u$ in $\Omega$. The subsolution $\underline u$ is a lower
bound of the solution. To obtain an upper bound, suppose $u$ attains
its maximum at $x_0\in \Omega$, that is $Du(x_0)=0$ and
$D^2u(x_0)\leqslant 0$. By ellipticity, we then have $A(x_0,0) < 0$,
which contradicts the hypothesis that $A(x,0)$ has a non-negative
eigenvalue for each $x\in\Omega$ and hence $u$ must take its maximum
on $\partial\Omega$. Thus we have the solution bound
\begin{equation}\label{4.1}
|u|\leqslant K_0,\hspace{2mm} {\rm{in}} \hspace{2mm} \bar\Omega,
\end{equation}
where $K_0$ depends on $\|\underline u\|_{L^\infty(\Omega)}$ and
$\|\varphi\|_{L^\infty(\Omega)}$.

The latter condition is equivalent to constant functions being
viscosity supersolutions, as defined in \cite{Tru1990}, so more
generally we could assume that there exists a viscosity
supersolution $\bar u \geqslant \phi$ on $\partial\Omega$.

We assume a quadratic bound from below for the matrix $A$ to control
the gradient of the solution, namely that $A$ satisfy the following
structure condition (\ref{S1}),
\begin{equation*}
A(x,p) \geqslant - \mu_0(1 + |p|^2)I.
\end{equation*}
Consider the function $\psi=e^{\kappa u}|Du|$ for some $\kappa>0$.
Suppose $\psi$ attains its maximum at $x_0\in \Omega$, namely $D_i
\psi(x_0)=0$. Hence, at the point $x_0$, we have $D_iuD_i \psi=0$.
By calculation, we have
\begin{equation}
\kappa |Du|^4+D_{ij}uD_iuD_ju=0.
\end{equation}
From the ellipticity condition $D^2u-A(x,Du)>0$, we have
\begin{equation}
\kappa |Du|^4+A_{ij}(x,Du)D_iuD_ju<0.
\end{equation}
By the structure condition (\ref{S1}) of the matrix $A$, we have
\begin{equation}\label{contradition ineqn}
\kappa |Du|^4-\mu_0(1+|Du|^2)|Du|^2<0.
\end{equation}
Without loss of generality, we may assume $Du(x_0)\geqslant 1$. The
left hand side of (\ref{contradition ineqn}) will be nonnegative if
we choose the positive constant $\kappa$ sufficiently large. This
contradiction leads to the gradient estimate
\begin{equation}\label{gradient bounds}
\sup\limits_\Omega |Du|\leqslant C,
\end{equation}
where $C$ depends on $\mu_0$, $\|u\|_{L^\infty(\Omega)}$ and
$\sup_{\partial \Omega}|Du|$.

On the boundary, the tangential derivatives of $u$ are given by the
boundary condition and the interior normal derivative bound from
below is controlled by the subsolution $\underline u$. The estimate
from above follows from condition (\ref{S1}) and the ellipticity,
which implies $\triangle u > A_{ii}$, (see proof of Theorem 14.1 in
\cite{GT2001}).
 Combining with
(\ref{gradient bounds}), we have
\begin{equation}\label{4.6}
|Du|\leqslant K_1,\hspace{2mm} {\rm{in}} \hspace{2mm} \bar\Omega,
\end{equation}
where $K_1$ depends on $\|\underline u\|_{C^1(\Omega)}$ and
$\|\varphi\|_{C^1(\Omega)}$.

As we shall indicate below, for the special cases of prescribed
Jacobian and optimal transportation equations, the above conditions
can be relaxed.

To complete the proof of Theorem 1.3, we have from Theorem 1.2 and
the $C^1$ bounds (\ref{4.1}) and (\ref{4.6}), uniform estimates in
$C^2(\bar\Omega)$ for classical elliptic solutions of the Dirichlet
problems:
\begin{equation}
\left\{ \begin{array}{rcll}
    \det\{D^2u-A(x,Du)\} & = & tB(x,Du) + (1-t)\det\{D^2\underline u-A(x,D\underline u)\}
     \ {\rm in} \  \Omega ,           \\
                       u & = & \varphi  \ {\rm on}\  \partial\Omega ,   \\
         \end{array}  \right.
\end{equation}
for $0 \leqslant t \leqslant 1$. Theorem 1.3 then follows from the
Evans - Krylov second derivative H\"older estimates of Evans, Krylov
and Caffarelli-Nirenberg-Spruck, (see for example \cite{GT2001},
Theorems 17.26, 17.26'), and the method of continuity
(\cite{GT2001}, Theorem 17.8). The uniqueness assertion is immediate
from the maximum principle.

\begin{remark}
We remark that the classical solution $u$ belongs to $C^\infty
(\bar\Omega)$ if both $A$ and $B$ are $C^\infty$ functions on
$\bar\Omega\times\mathbb{R}^n$, as well as $\partial\Omega\in
C^\infty, \varphi\in C^\infty(\partial\Omega)$. Under the stated
hypotheses it also follows from the linear theory \cite{GT2001} that
the solution $u\in C^{3,\alpha}(\bar\Omega)$ for all $\alpha<1$. In
fact we have the stronger inclusion, $u\in W^{4,p}(\Omega)$ for all
$p<\infty$ and this would in fact suffice for the argument in
Section 2.
 Note that we only use $u\in C^3(\bar\Omega)$ in Section 3.
\end{remark}

\begin{remark}
As a special case for Monge-Amp\`ere type equation (\ref{Dirichlet
BVP of OTE}), the existence results for the classical solutions to
the Dirichlet problem for the equation $\det
(u_{ij}+\sigma_{ij}(x))=\psi(x)$ treated in the strictly convex
domain in \cite{CNSI84,CNS87} can be generalized to non-convex
domains. Also, the equation with general right hand side treated by
Li \cite{Liyanyan90} can be generalized to non-convex domains.
\end{remark}

\begin{remark} For the relationship with equations arising in conformal geometry
the reader is referred to \cite{Guanbo2007} where a similar barrier
argument for boundary estimates is used.
\end{remark}

In the remainder of this paper, we discuss the application to
optimal transportation and prescribed Jacobian equations. Optimal
transportation problems have received a lot of attention in recent
years and for the basic theory, we refer to the the books,
\cite{Villani03,Villani08}.
 Letting $Y$ be a $C^1$ mapping from
$\Omega\times\mathbb{R}\times\mathbb{R}^n$ to $\mathbb{R}^n$ and $u$
a function in $C^2(\Omega)$, we define a mapping $T =
Tu:\Omega\rightarrow \mathbb{R}^n$ by setting $Tu=Y(\cdot,u,Du)$.
The corresponding prescribed Jacobian equation may be written  the
form
\begin{equation}\label{4.8}
|\det DT|= \psi(\cdot,u,Du)
\end{equation}
for a given positive function $\psi$ on
$\Omega\times\mathbb{R}\times\mathbb{R}^n$. If $\det D_pY\ne 0$, we
can then write equation (\ref{4.8}), for elliptic solutions $u$, as
an equation of Monge-Amp\`ere type
\begin{equation}\label{4.9}
\det [D^2u-A(\cdot,u,Du)]=B(\cdot,u,Du)
\end{equation}
with matrix function $A$ and scalar function $B$ given by
\begin{equation}\label{4.10}
\begin{array}{lll}
A(\cdot,u,p)&=&-Y^{-1}_p(Y_x+Y_u\otimes p),\\
B(\cdot,u,p)&=&|\det Y_p|^{-1}\psi,
\end{array}
\end{equation}
so that when $Y$ and $\psi$ are independent of $u$, (and $B$ is
positive), we obtain a Monge-Amp\`ere type equation of the form
(\ref{Dirichlet BVP of OTE}) considered here. Optimal transportation
equations are the special cases where the mapping $Y$ is generated
by a cost function $c$ satisfying conditions (A1) and (A2) as
formulated in the introduction. We then have, in accordance with
(\ref{OT}),
\begin{equation}\label{4.11}
\begin{array}{rll}
Y^{-1}_p  &=& D^2_{x,y}c(\cdot,Y),\\
A(\cdot,p)&=&D^2_xc(\cdot,Y),\\
B(\cdot,p)&=&|\det D^2_{x,y}c|\psi >0.
\end{array}
\end{equation}

The natural boundary condition for prescribed Jacobian equations is
the prescription of the image $T(\Omega)$. The existence of
classical solutions is treated in \cite{Tru2008}, which extends the
optimal transportation case in \cite{TruWang2009}. However the
classical Dirichlet problem in small balls has been used for local
regularity arguments in \cite{MTW2005,TruWang2009C1,Wang1996}.
Clearly, the estimates in Theorem \ref{Th1} and \ref{Th2} extend
immediately to embrace prescribed Jacobian and optimal
transportation equations under regularity of the matrix functions
$A$ given by (\ref{4.10}) and (\ref{4.11}). For the above solution
and gradient estimates some relaxation of our additional hypotheses
is possible. First we note that for any vector $\xi\in
\mathbb{R}^n$, the quantity $Tu\cdot \xi$ must assume its maximum
and minimum values on the boundary $\partial \Omega$ since $\det DT
\ne 0$. In the optimal transportation case this immediately provides
a bound for the gradient in term of its boundary trace since
$Du=c_x(\cdot, Tu)$ in $\Omega$. Also in the optimal transportation
case, functions of the form $\bar u=c(\cdot, y_0)$ are solutions of
the homogeneous equation so we automatically obtain solution bounds.
For an arbitrary domain $\Omega$, a condition of the type (\ref{S1})
would still be required for the boundary gradient estimate or more
generally the existence of a supersolution $\bar u$ satisfying $\bar
u=\varphi$ on $\partial \Omega$. However if the solution is
$c$-convex in $\Omega$, that is at each $x_0\in \Omega$, there
exists $y_0\in \mathbb{R}^n$, such that
\begin{equation}\label{4.12}
u(x)\geqslant u(x_0)+c(x,y_0)-c(x_0,y_0)
\end{equation}
in $\Omega$, then we have $Tu(\Omega)\subset T\underline u (\Omega)$
by virtual of the monotonicity property (\cite{MTW2005}, Lemma 4.3),
and the global gradient bound follows immediately. From
\cite{TruWang2009C1}, we know that an elliptic solution will be
$c$-convex in $\Omega$ if $\Omega$ is $c$-convex with respect to
each $y\in \mathbb{R}^n$, that is the images $c_y(\cdot, y)(\Omega)$
are convex in $\mathbb{R}^n$. Consequently, we have the following
existence theorem for optimal transportation equations.
\begin{theorem}\label{Th4.1}
Let $c\in C^4(\mathbb{R}^n\times \mathbb{R}^n)$ be a cost function
satisfying (A1), (A2) with regular matrix function $A$ given by
(\ref{OT}) and $B\in C^2(\bar \Omega\times\mathbb{R}^n)$ satisfying
$\inf B>0$. Suppose there exists a subsolution $\underline u\in
C^2(\bar \Omega)$ of equation (\ref{Dirichlet BVP of OTE}) in
$\Omega$, satisfying $\underline u=\varphi$ on $\partial \Omega$.
Then if there also exists a supersolution $\bar u$, $=\varphi$ on
$\partial \Omega$, or if $\Omega$ is $c$-convex with respect to all
$y\in \mathbb{R}^n$, there exists a unique classical solution $u\in
C^2(\bar \Omega)$ of the Dirichlet problem (\ref{Dirichlet BVP of
OTE}) with $u\geqslant \underline u$ in $\Omega$.
\end{theorem}

\begin{remark} A concept of generalized solution for optimal transportation
equations, extending that of Aleksandrov for the standard
Monge-Amp\`ere equation, is also introduced in \cite{MTW2005}  and
the Dirichlet problem for generalized solutions is considered in
\cite{GN2007}. Without the regularity assumption, this notion can be
non-local as shown by Loeper \cite{Loeper09Acta} and thus is not
really a weak form of the classical notion.
\end{remark}

\begin{remark}
These results also extend to the more general Monge-Amp\`ere type
equations (\ref{4.9}) for $A$ and $B$ satisfying
\begin{equation*}
D_uA_{ij}\xi_i\xi_j\geqslant 0, \  B_u\geqslant 0
\end{equation*}
for all $(x,u,p) \in \Omega\times \mathbb{R}\times \mathbb{R}^n$. In
this generality, Theorem \ref{Th4.1} also extends to embrace
prescribed Jacobian equations with mappings $Y$ determined by
generating functions as in \cite{Tru2012}.
\end{remark}


\end{document}